\documentclass[10pt]{article}

\usepackage {amssymb}
\newcommand{\n}{\noindent}
\renewcommand{\a}{\alpha}

\renewcommand{\b}{\beta}

\renewcommand{\o}{\omega}

\newcommand{\f}{\mathfrak}

\begin{document}

\section*{Description de la structure de certaines
superalg\`ebres de Lie quadratiques via la notion de la $T^*-$extension}

\vspace{0.5cm}

Ignacio BAJO,~Sa\"{\i}d BENAYADI~et~Martin BORDEMANN

\vspace{0.4cm}

\noindent I.B.: Universidad de Vigo, Depart. Matem\'atica Aplicada, \\
E.T.S.I Telecomunicaci\'on, 36280 Vigo, Espagne.\\
E-mail: ibajo@dma.uvigo.es

\vspace{0.1cm}

\noindent S.B.: Universit\'e de Metz, D\'epartement de Math\'ematiques,\\
CNRS, UPRES-A-7035,
Ile du Saulcy, F-57045 Metz cedex 1, France.\\
E-mail: benayadi@poncelet.sciences.univ-metz.fr

\vspace{0.1cm}

\noindent M.B.: Fakult{\"a}t f{\"u}r Physik, Universit{\"a}t Freiburg,\\
Hermann-Herder-Str.3, D-79104 Freiburg i. Br., Allemagne.\\
E-mail: Martin.Bordemann@physik.uni-freiburg.de

\vspace{1cm}

\n {\bf R\'esum\'e.} Dans cette note, nous introduisons la notion de
$T^*-$extension $T^*{\f g}$ d'une superalg\`ebre de Lie ${\f g}$, 
c'est-\`a-dire une
extension de ${\f g}$ par son espace dual ${\f g}^*$. L'accouplement 
naturel induit
sur cette extension une forme bilin\'eaire, paire, supersym\'etrique, non
d\'eg\'en\'er\'ee et invariante
($B([X,Y],Z)=B(X,[Y,Z])$ pour tous $X,Y,Z \in T^*{\f g}$), c'est-\`a-dire
la structure
d'une superalg\`ebre de Lie quadratique (ou orthogonale). Ces extensions
de ${\f g}$ se classifient
\`a l'aide de son troisi\`eme groupe de cohomologie scalaire paire. En
outre, nous
montrons que toutes les superalg\`ebres de Lie quadratiques
${\f a}={\f a}_{\bar{0}} \oplus {\f a}_{\bar{1}}$ de dimension finie qui
sont soit nilpotentes,
soit r\'esolubles telles que
$[{\f a}_{\bar{1}},{\f a}_{\bar{1}}]\subset [{\f a}_{\bar{0}},{\f a}_{\bar 0}]$
s'obtiennent \`a l'aide d'une $T^*-$extension dans le cas d'un corps
alg\'ebriquement clos de caract\'eristique nulle.

\subsection*{Description of the structure of quadratic Lie
superalgebras via the notion of $T^*-$extension}

{\bf Abstract.} In this note we introduce the notion of $T^*-$extension
$T^*{\f g}$ of a Lie superalgebra ${\f g}$, i.e.~an extension of ${\f g}$ 
by its
dual space ${\f g}^*$. The natural pairing induces on $T^*{\f g}$ an
even supersymmetric
nondegenerate bilinear form $B$ which is invariant
($B([X,Y],Z)=B(X,[Y,Z])$ for all $X,Y,Z \in T^*{\f g}$), i.e.~the
structure of a quadratic
(or metrised or orthogonal) Lie superalgebra. These extensions can be
classified by the
third even scalar cohomology group of ${\f g}$. Moreover, we show that
all finite-dimensional
quadratic Lie superalgebras ${\f a}={\f a}_{\bar{0}} \oplus {\f a}_{\bar{1}}$
which are either nilpotent, or solvable and such that
$[{\f a}_{\bar{1}},{\f a}_{\bar{1}}]\subset
[{\f a}_{\bar{0}},{\f a}_{\bar{0}}]$ can be constructed
by means of a $T^*-$extension in the case of an algebraically closed
field of characteristic zero.

\subsection*{Abridged English Version.}

All the Lie superalgebras in this note will be finite-dimensional over an
algebraically field $\mathbb K$ of characteristic zero. For the
terminology used here we refer the reader to [4] and [6].

For a given Lie superalgebra ${\f g}= {\f g}_{\bar{0}} \oplus {\f g}_{\bar{1}}$
a bilinear form $B$ on $\f g$ is called an invariant scalar product if and
only if $B$ is even, supersymmetric, nondegenerate, and invariant, the
latter meaning that $B([X,Y],Z)= B(X,[Y,Z])$ for all $X, Y, Z\in{\f g}$.
The pair $(\f g, B)$ is called a quadratic (or metrised or orthogonal) Lie 
superalgebra in
case $B$ is an invariant scalar product on $\f g$.
$\varphi:(\f g, B)\rightarrow (\f g', B')$ is called an isometry of
the quadratic Lie superalgebras $(\f g, B)$ and $(\f g', B')$ if and only if
$\varphi$ is an isomorphism of the Lie superalgebras $\f g$ and $\f g'$
and $B'(\varphi X,\varphi Y)=B(X,Y)$ for all $X,Y\in \f g$.

In [1] Benamor and Benayadi have transferred the notion of double extension
introduced for quadratic Lie algebras by Medina and Revoy [5] to certain classes
of Lie superalgebras (those where the centre has nontrivial intersection 
with the
even part) which allows an inductive construction of those quadratic Lie
superalgebras. The aim of this note is to transfer the notion of $T^*$-extension
formulated in [3] for quadratic Lie algebras to the case of a quadratic
Lie superalgebra.

Let ${\f g}= {\f g}_{\bar{0}} \oplus {\f g}_{\bar{1}}$ be a Lie
superalgebra, $[~,~]_{\f g}$ its bracket,  ${\f g}^*$ its
dual space, and $\pi$ its coadjoint representation. Let $\o$ be an even
$2$-cocycle of ${\f g}$ with values in ${\f g}^*$.

\n {\bf Proposition I.1./I.2.} With the above notation define the following
structures on the vector space ${\mathcal A}:= {\f g}\oplus {\f g}^*$:
Let $~X+F$ (resp. $Y+H)$ be a homogeneous element of $\mathcal A$ of
degree $x$ (resp. $y$) in $\mathbb Z /2\mathbb Z$, and
\begin{eqnarray*}
    [X+F,Y+H]   & := & [X,Y]_{\f g}+\o(X,Y)+\pi(X)(H)-(-1)^{xy}\pi(Y)(F), \\
     B(X+F,Y+H) & := &  F(Y)~+~(-1)^{xy}~H(X).
\end{eqnarray*}
Then $({\mathcal A},B)$ is a quadratic Lie superalgebra if and only if
$\o$ is supercyclic, i.e.
\[
  \o(X,Y)(Z)=~ (-1)^{x(y+z)}~\o(Y,Z)(X),~~\forall~(X,Y,Z) \in
                                    {{\f g}_x\times {\f g}_y \times {\f g}_z}.
\]

We shall speak of the quadratic Lie superalgebra $({\mathcal A},B)$
constructed out of $\f g$ and the supercyclic $2$-cocycle $\o$ as a
$T^*$-extension of $\f g$ and shall denote ${\mathcal A}$ by the symbol
$T^*_\o {\f g}$ or $T^* {\f g}$.

\n {\bf Theorem I.2.} Let ${\f g}$ a Lie superalgebra and $\o$ be a
supercyclic $2$-cocycle of ${\f g}$ with values in the dual space
${\f g}^*$.\\
Then $\hat{\o}:{\f g}\times {\f g}\times {\f g}\rightarrow {\mathbb K}$
defined by $\hat{\o}(X,Y,Z):=\o(X,Y)(Z)$ for all $X,Y,Z\in {\f g}$ is an
even
$3$-cocycle of ${\f g}$ with values in the trivial ${\f g}$-module
${\mathbb K}$, and each even scalar $3$-cocycle of ${\f g}$ is given by
some $\hat{\o}$.
Moreover two $T^*$-extensions $T^*_{\o_1} {\f g}$ and $T^*_{\o_2} {\f g}$
are isometric if the corresponding even scalar $3$-cocycles $\hat{\o}_1$ and
$\hat{\o}_2$ are cohomologous.

The main result of this note is the following:

\n {\bf Theorem II.1.} Let $({\f g},B)$ be a quadratic Lie superalgebra which
is either nilpotent, or solvable and such that
$[{\f g}_{\bar{1}},{\f g}_{\bar{1}}] 
\subset [{\f g}_{\bar{0}},{\f g}_{\bar{0}}]$. \\
Then $\f g$
contains a graded totally isotropic ideal ${\f I}$
(i.e. $B({\f I},{\f I})=\{0\}$) whose dimension is equal to the
integer part of one half of the
dimension of ${\f g}$.\\
In case the dimension of $\f g$ is even then
$\f g$ is isometric to the $T^*-$extension of the quotient Lie
superalgebra ${\f g}/{\f I}$. \\
In case the dimension of $\f g$ is odd then
$\f g$ is isometric to a graded ideal of codimension $1$ of the
$T^*-$extension of the quotient Lie superalgebra ${\f g}/{\f I}$
restricted to which the invariant scalar product is nondegenerate.

\n {\bf Remarks:} 1. The proof of this Theorem relies on the Theorems of
Engel and Lie where the latter is known to no longer hold in the original form
for all solvable Lie superalgebras but only for the class mentioned above.

2. There is an example of a nilpotent quadratic Lie superalgebra whose
centre is entirely contained in its odd part; it is a $T^*$-extension, but
does not fall under the class treated by [1].

\subsection*{Introduction}

Les superalg\`ebres de Lie envisag\'ees dans cette note
sont de
dimension finie sur un corps  $\mathbb K$ commutatif alg\'ebriquement clos de
caract\'eristique z\'ero.

Soit ${\f g}= {\f g}_{\bar{0}} \oplus {\f g}_{\bar{1}}~$  une superalg\`ebre
de Lie, une forme bilin\'eaire $B$ sur $\f g$ est dite invariante si:
$B([X,Y],Z)= B(X,[Y,Z])$ pour tous $X, Y, Z$ \'el\'ements de $\f g$, si 
en plus $B$ est
supersym\'etrique (c'est-\`a-dire $~B(X,Y)= (-1)^{xy} B(Y,X)$ $\forall~(X,Y)
\in
{\f g}_x \times {\f g}_y),~$ paire  (c'est-\`a-dire
$~B({\f g}_{\bar{0}},{\f g}_{\bar{1}})= \{0\})~$   et non d\'eg\'en\'er\'ee,
$B$ est alors dite un produit scalaire
invariant sur $\f g$. Une superalg\`ebre de Lie $~{\f g}~$ est dite
quadratique si
elle est munie d'un produit scalaire invariant. Un id\'eal gradu\'e $\f I$
d'une
superalg\`ebre de Lie quadratique $~({\f g},B)~$ est dit non d\'eg\'en\'er\'e
(resp.
d\'eg\'en\'er\'e) si la restriction de $~B~$ \`a $~{\f I}\times {\f I}~$ est non
d\'eg\'en\'er\'ee (resp. d\'eg\'en\'er\'ee). Une superalg\`ebre de Lie 
quadratique
$~({\f g},B)~$ est dite $~B$-irr\'eductible si $~{\f g}$ ne contient aucun
id\'eal gradu\'e non d\'eg\'en\'er\'e diff\'erent de $~\{0\}~$ et de  
$~{\f g}$.  Deux
superalg\`ebres de Lie quadratiques  $~({\f g},B)~$ et $~({\f g'},B')~$
sont isom\'etriques s'il existe $~\varphi~$ un isomorphisme de 
superalg\`ebres de Lie de
$~{\f g}~$ sur  $~{\f g'}~$ tel que  $~B'(\varphi(X),\varphi(Y))=
B(X,Y),~~\forall~(X,Y) \in  {\f g}\times {\f g}.~$

Dans $[1]$, il est montr\'e que toute superalg\`ebre de Lie quadratique
$~({\f g},B)~$   $~B$-irr\'eductible de dimension $n$ telle que
la partie paire du centre ${\f z}(\f g)$ soit non nulle
est une double extension d'une
superalg\`ebre de Lie quadratique de dimension $n-2$ par l'alg\`ebre de Lie
de dimension $1$. En utilisant ce r\'esultat, les auteurs de $[1]$ ont
obtenu une classification inductive des superalg\`ebres de Lie
quadratiques $~({\f g}= {\f g}_{\bar{0}} \oplus {\f g}_{\bar{1}},B)~$
telles que dim${\f g}_{\bar{1}}= 2~$ en montrant que de telles
superalg\`ebres v\'erifient $~{\f z}(\f
g)\cap {\f g}_{\bar{0}}\neq~ \{0\}.~$ Rappelons que la notion de la
double extension a \'et\'e introduite dans $[5]$  par A. M\'edina et Ph.
Revoy pour \'etudier les alg\`ebres de Lie quadratiques et a \'et\'e
g\'en\'eralis\'ee par H.Benamor et S. Benayadi aux superalg\`ebres de Lie.
Gr\^ace \`a cette notion,  une
classification inductive de toutes les alg\`ebres de Lie
quadratiques est obtenue dans $[5]$ et une classification
 inductive de certaines superalg\`ebres de Lie quadratiques est
obtenue dans $[1]$ et dans $[2]$.

Si $~({\f g},B)~$ est une superalg\`ebre de Lie quadratique telle que
${\f z}(\f g)\neq \{0\}~$ et
$~{\f z}(\f g)\cap {\f g}_{\bar{0}}=~ \{0\}~$ c'est-\`a-dire
$~{\f z}(\f g)\subset {\f g}_{\bar{1}},~$ les arguments utilis\'es dans $[1]$
ne marchent
plus pour montrer que $\f g$ est une double extension. Par
cons\'equent, il faut chercher d'autres m\'ethodes pour d\'ecrire au moins
des sous-classes de la classe $\mathcal C$ des superalg\`ebres de Lie
quadratiques $~({\f g},B)~$ telles que ${\f z}(\f g)\neq \{0\}~$
et  $~{\f z}(\f g)\subset {\f g}_{\bar{1}}.~$
Dans le premier paragraphe de cette note, nous
allons montrer que la classe $\mathcal C$ est non vide en donnant des
exemples d'\'el\'ements de $\mathcal C$ qui sont nilpotents.

Dans cette note, nous allons g\'en\'eraliser la notion de la $T^*-$extension
aux
superalg\`ebres de Lie, cette notion a \'et\'e introduite par M. Bordemann
dans
$[3]$ pour \'etudier les alg\`ebres non associatives munies de formes 
bilin\'eaires, sym\'etriques, non d\'eg\'en\'er\'ees et associatives (ou 
invariantes), donc  en
particulier pour  \'etudier les alg\`ebres de Lie quadratiques. Apr\`es, nous
allons
utiliser cette notion pour d\'ecrire les superalg\`ebres de Lie quadratiques
nilpotentes et les superalg\`ebres de Lie r\'esolubles
$~({\f g}= {\f g}_{\bar{0}} \oplus {\f g}_{\bar{1}},B)~$
telles que
$~[{\f g}_{\bar{1}},{\f g}_{\bar{1}}]\subset 
[{\f g}_{\bar{0}},{\f g}_{\bar{0}}]~$.
Rappelons que si $~({\f g},B)~$ est
une superalg\`ebre de Lie quadratique alors le sous-espace orthogonal de
$~[\f g,\f g]~$
par rapport \`a $B$ est ${\f z}(\f g).~$ Par cons\'equent, si
$~({\f g},B)~$ est une superalg\`ebre de Lie quadratique r\'esoluble
diff\'erente de $\{0\}$ alors $~{\f z}(\f g)\neq \{0\}.~$

\subsection*{I. $T^*$-extensions des superalg\`ebres de Lie quadratiques.}

Soient $~{\f g}= {\f g}_{\bar{0}} \oplus {\f g}_{\bar{1}}~$   une
superalg\`ebre de Lie, $~[~,~]_{\f g}~$ son crochet,   $~{\f g}^*~$ son
dual et
$~\pi~$ sa repr\'esentation coadjointe. Rappelons que $\pi$ est d\'efinie
par: $(\pi(X)(F))(Y)= -(-1)^{xf} F([X,Y])~$ o\`u $~X ~$(resp. $Y)~$ est un
\'el\'ement homog\`ene
de $\f g$ de degr\'e $x ~$(resp. $y)$ et $~F~$ est un \'el\'ement homog\`ene de
${\f g}^*~$
de degr\'e $f$. Soit $\o: {\f g}\times {\f g} \rightarrow {\f g}^*~$ une
application bilin\'eaire paire, c'est-\`a-dire
$~\o({\f g}_x \times{\f g}_y)\subset {\f g}^*_{x+y}~~
\forall~(x,y)\in {\mathbb Z}/2{\mathbb Z}\times {\mathbb Z}/2{\mathbb Z}$.

\noindent On d\'efinit sur l'espace vectoriel
${\mathcal A}:= {\f g}\oplus {\f g}^*~$
l'application bilin\'eaire
$~[~,~]: {\mathcal A}\times {\mathcal A} \rightarrow {\mathcal A}~$
d\'efinie par:
\[
[X+F,Y+H]=~[X,Y]_{\f g}~+~\o(X,Y)~+~\pi(X)(H)~-~(-1)^{xy}\pi(Y)(F),
\]

\n o\`u $~X+F$ (resp. $Y+H)~$ est un \'el\'ement homog\`ene de $~\mathcal A~$ de
degr\'e $x$ (resp. $y).~$
Cette multiplication fait de $~\mathcal A~$ une alg\`ebre
${{\mathbb Z}/ {2{\mathbb Z}}}-$gradu\'ee.

\vspace{0.2cm}

\n {\bf Proposition I.1.} $~({\mathcal A},[~,~])~$ est une superalg\`ebre de Lie
si et
seulement si $~\o ~\in~(Z^2({\f g},{\f g}^*))_{\bar{0}},~$ c'est-\`a-dire
$~\o~$ est
paire, superantisym\'etrique ($\omega(X,Y)=(-1)^{xy}\omega(Y,X)$ pour tous
$X,Y\in {\f g}_{x}\times {\f g}_{y}$) et v\'erifie
\begin{eqnarray*}
   \o(X,[Y,Z])+(-1)^{x(y+z)}\o(Y,[Z,X])
   +(-1)^{z(x+y)}\o(Z,[X,Y])
                 +\pi(X)(\o (Y,Z))         &   & \\
 +(-1)^{x(y+z)}\pi(Y)(\o (Z,X))
   +(-1)^{z(x+y)}\pi(Z)(\o (X,Y))
                                           & = & ~0,
\end{eqnarray*}
pour tous $(X,Y,Z) \in {{\f g}_x\times {\f g}_y \times {\f g}_z}$.
Sur ${\mathcal A}= {\f g}\oplus {\f g}^*,~$ on consid\`ere la forme bilin\'eaire
$~B:{\mathcal A}\times {\mathcal A} \rightarrow {\mathbb K}~$ d\'efinie par:
\[
B(X+F,Y+H):=~ F(Y)~+~(-1)^{xy}~H(X),~~\forall~ (X+F,Y+H)~\in
                                         ~{\mathcal A}_x \times{\mathcal A}_y.
\]
Cette forme est supersym\'etrique, paire et non d\'eg\'en\'er\'ee.

\vspace{0.2cm}

\n {\bf Proposition I.2.} $B~$ est invariante si et seulement si $~\o~$
est supercyclique,  c'est-\`a-dire $~\o~$ v\'erifie
\[
\o(X,Y)(Z)=~ (-1)^{x(y+z)}~\o(Y,Z)(X),~~\forall~(X,Y,Z) \in
                                    {{\f g}_x\times {\f g}_y \times {\f g}_z}.
\]

\vspace{0.2cm}

\n {\bf D\'efinition.} Soient $\f g$ une superalg\`ebre de Lie et
$~\o ~\in~(Z^2({\f g},{\f g}^*))_{\bar{0}}~$ qui est en plus supercyclique. La
superalg\`ebre de Lie  $~{\mathcal A}= {\f g}\oplus {\f g}^*~$ munie du
crochet
$~[~,~]~$ et de la forme bilin\'eaire $~B~$ d\'efinis ci-dessus est appel\'ee
la
$T^*-$extension de $~\f g~$ par $~\o~$, et on la note $T^*_{\o}{\f g}$
ou
$~(T^*_{\o}{\f g},B)$ ou parfois $T^*{\f g}$.

\vspace{0.2cm}

\n {\bf Remarque.} Si $\o= 0$, $T^*_{0}{\f g}$ n'est autre que la double
extension de $\{0\}$ par ${\f g}.~$

\vspace{0.2cm}

\n {\bf Lemme I.1.} Soit $({\f g},B)$ une superalg\`ebre de Lie
quadratique de dimension paire $n$.

\n Alors, un sous-espace vectoriel $\f I$ totalement isotrope de $\f g$ de
dimension $n\over 2$ est un id\'eal de $\f g$ si et seulement si $\f
I$ est ab\'elien,  c'est-\`a-dire $[\f I,\f I]= \{0\}$.

\vspace{0.2cm}

Le th\'eor\`eme suivant est le premier r\'esultat principal de cette note.

\vspace{0.1cm}

\n {\bf Th\'eor\`eme I.1.}  Soit $({\mathcal A},T)$ une superalg\`ebre de Lie
quadratique de dimension $n$.

\n Alors, $({\mathcal A},T)$ est isomorphe \`a une $T^*-$extension
$(T^*_{\o}{\f g},B)$ si et seulement si $n$ est pair et $\mathcal A$ contient un
id\'eal gradu\'e $\f I$ totalement isotrope de dimension $n\over 2$.
Dans ce cas $\f g$ et ${\mathcal A}/{\f I}$ sont isomorphes en tant
que superalg\`ebres de Lie.

\vspace{0.2cm}

Nous avons la d\'escription cohomologique suivante des $T^*-$extensions:

\vspace{0.1cm}

\n {\bf Th\'eor\`eme I.2.}  Soit $\f g$ une superalg\`ebre de Lie sur un corps
${\mathbb K}$.

\n i) Pour tout $~\o ~\in~(Z^2({\f g},{\f g}^*))_{\bar{0}},~$
l'application canonique $\o~\mapsto~{\hat \o}:~((X,Y,Z)~\mapsto~\o(X,Y)(Z)), $
o\`u $(X,Y,Z) ~\in~{\f g}\times {\f g}\times {\f g},~$ d\'efinit un
isomorphisme lin\'eaire entre le
sous-espace de tous les \'el\'ements supercycliques de
$(Z^2({\f g},{\f g}^*))_{\bar{0}}$ et l'espace
 $(Z^3({\f g},{\mathbb K}))_{\bar{0}}$ de tous les $3-$cocycles scalaires
paires de $\f g$. Rappelons
que $(Z^3({\f g},{\mathbb K}))_{\bar{0}}$ est l'espace vectoriel des formes
trilin\'eaires $f: {\f g}\times {\f g}\times {\f g} \rightarrow {\mathbb K}$
qui sont paires ($f({\f g}_x,{\f g}_y,{\f g}_z)=\{0\}$ si $x+y+z=0$),
superantisym\'etriques
($f(X,Y,Z)=-(-1)^{xy}f(Y,X,Z)=-(-1)^{yz}f(X,Z,Y)$ pour tous
$(X,Y,Z)\in {\f g}_x\times {\f g}_y \times {\f g}_z$) et ferm\'ees:
\begin{eqnarray*}
           0 & =  & f([X,Y],Z,V) -(-1)^{yz}f([X,Z],Y,V) + 
                                         (-1)^{x(y+z)}f([Y,Z],X,V) \\
             &    &    +(-1)^{(y+z)v}f([X,V],Y,Z) 
                                         -(-1)^{x(y+v)+vz}f([Y,V],X,Z) \\
             &    &           + (-1)^{(x+y)(z+v)}f([Z,V],X,Y) ,
\end{eqnarray*}
pour tous
$(X,Y,Z,V)\in {\f g}_x\times {\f g}_y \times {\f g}_z \times{\f g}_v$.

\n ii) Soit $\phi: {\f g} \times {\f g} \rightarrow {\mathbb K}$ une
application bilin\'eaire paire ($\phi({\f g}_x,{\f g}_y)=\{0\}$ si 
$x+y\neq 0$) et
superantisym\'etrique, et soient $\o_1$ et $\o_2$ deux \'el\'ements 
supercycliques
de $(Z^2({\f g},{\f g}^*))_{\bar{0}}$.
  Alors l'application
\[
 S_\phi : T^*_{\o_1}{\f g} \rightarrow T^*_{\o_2}{\f g}:
               (X+F)\mapsto X+ \phi(X,.) + F
\]
d\'efinit une isom\'etrie de superalg\`ebres de Lie quadratiques si et
seulement si $\o_2=\o_1-\delta\phi$ o\`u
$(\delta\phi)(X,Y,Z):= -\phi([X,Y],Z) +(-1)^{yz}\phi([X,Z],Y)
-(-1)^{x(y+z)}\phi([Y,Z],X)$
pour tous $(X,Y,Z)\in {\f g}_x\times {\f g}_y \times {\f g}_z $.

\vspace{0.1cm}

\n Ainsi on obtient une application naturelle
$[\hat{\omega}]\mapsto [T^*_\o {\f g}]$ du troisi\`eme groupe de cohomologie
scalaire paire de $\f g$ dans l'ensemble des classes d'isom\'etries des 
structures
de superalg\`ebres de Lie quadratiques de dimension $2\dim \f g$.

\vspace{0.2cm}
\vspace{0.1cm}

\n {\bf Exemple I.1.} Maintenant nous allons construire des superalg\`ebres
de Lie quadratiques $(\f g,B)$ qui sont nilpotentes et qui v\'erifient
${\f g}_{\bar{0}}\neq\{0\}$ et
${\f g}\in {\mathcal C}$.
Avant de construire ces \'el\'ements de $\mathcal C$,
remarquons que $\mathcal C$ contient toutes les superalg\`ebres de Lie $\f g$
telles que ${\f g}_{\bar{0}}= \{0\}$.

\n Soit ${\mathcal T}(n)$ (resp.~${\mathcal N}(n)$) le sous-espace vectoriel de
${\mathcal M}_n(\mathbb K)$, o\`u $n \in {\mathbb N}$,  form\'e des matrices
carr\'ees d'ordre
$n$ qui sont triangulaires sup\'erieures (resp.~triangulaires sup\'erieures
strictes). Il est facile de v\'erifier que le ${\mathbb K}-$espace vectoriel
\[
 {\f g}(n):= \{ \left(\begin{array}{cc}
                        A  &   B\\
                        C  &   D
                                 \end{array}\right)~:~(A,C,D) \in
                {{\mathcal N}(n)}^3, B\in {\mathcal T}(n)\}
\]
est une
sous-superalg\`ebre de Lie nilpotente de ${\f g}{\f l}(n,n)$ telle que
dim$~{\f z}({\f g}(n))=~ 1,~$
${\f z}({\f g}(n))\subset ({\f g}(n))_{\bar{1}},~$ et
$~[({\f g}(n))_{\bar{1}},({\f g}(n))_{\bar{1}}]=({\f g}(n))_{\bar{0}}$.
Rappelons que le crochet de ${\f g}{\f l}(n,n)$
est d\'efini
par $[M,N]= MN-(-1)^{\a\b}NM$, pour tout
$(M,N)\in {{\f g}{\f l}(n,n)}_{\a}\times  {{\f g}{\f l}(n,n)}_{\b}$.
Consid\'erons maintenant
${\mathcal E}=T^*_0({\f g}(n))$ la $T^*-$extension de ${\f g}(n)$ par $\o= 0$.
Le centre ${\f z}({\mathcal E})$ de cette $T^*-$extension est \'egale la somme
de ${\f z}({\f g}(n))$ et $ \{f \in ({\f g}(n))^*:
         ~f([{\f g}(n),{\f g}(n)])= \{0\}\}$.
Comme $~[({\f g}(n))_{\bar{1}},({\f g}(n))_{\bar{1}}]=({\f g}(n))_{\bar{0}},$
alors
${\f z}({\mathcal E}) \subset
{\f z}({\f g}(n))\oplus (({\f g}(n)^*)_{\bar{1}}= \{f \in ({\f g}(n))^*:
   ~f(({\f g}(n))_{\bar{0}})= \{0\}\}),~$
par cons\'equent
${\f z}({\mathcal E}) \subset
     ({\f g}(n))_{\bar{1}}\oplus ({\f g}(n)^*)_{\bar{1}}= 
                   {\mathcal E}_{\bar{1}},~$
en plus
${\f z}({\mathcal E})\neq \{0\}$ car ${\f z}({\f g}(n))\neq \{0\}.$
En utilisant la d\'efinition du crochet de $\mathcal E$, un calcul 
simple montre que
la nilpotence de ${\f g}(n)$ entra\^{\i}ne la nilpotence de $\mathcal E$. 
On conclut
que $\mathcal E$ est une superalg\`ebre de Lie quadratique telle que
${\f z}({\mathcal E})\neq\{0\}$ et ${\mathcal E}\in {\mathcal C}$.

\vspace{0.2cm}
\vspace{0.2cm}

\subsection*{II. Structure de certaines superalg\`ebres de Lie quadratiques
            r\'esolubles.}

\n {\bf Lemme II.1.} Soient $V$ un espace vectoriel
${{\mathbb Z}/{2{\mathbb Z}}}-$gradu\'e
muni d'une forme bilin\'eaire $B$ supersym\'etrique paire non 
d\'eg\'en\'er\'ee et
$\mathcal L$ une sous-superalg\`ebre de Lie de ${\f  o}{\f s}{\f p}(V,B)$. Si
$\mathcal L$ est constitu\'ee d'endomorphismes nilpotents de $V$ ou si
$\mathcal L$ est r\'esoluble telle que
$[{\mathcal L}_{\bar{1}},{\mathcal L}_{\bar{1}}]\subset
                [{\mathcal L}_{\bar{0}},{\mathcal L}_{\bar{0}}]$,
alors tout sous-espace vectoriel gradu\'e $W$ de $V$ totalement
isotrope et stable par $\mathcal L$ est contenu dans un sous-espace
vectoriel gradu\'e
$W_{max}$ de $V$ totalement isotrope maximal parmi les sous-espaces
vectoriels totalement isotropes de $V$ qui sont stables par $\mathcal L$
et  dim$~W_{max}= E({n\over 2})$ o\`u $E({n\over 2})$ est la partie
enti\`ere de
$n\over 2$. Si $n$ est pair, $W_{max}= (W_{max})^\perp$. Si $n$ est
impair, $W_{max}\subset (W_{max})^\perp$ 
et dim$~(W_{max})^\perp~-~$dim$~W_{max} =~1.$ On a aussi
$\varphi ((W_{max})^\perp)$ est contenu dans $W_{max}$, pour tout
$\varphi$ \'el\'ement de $\mathcal L$.

\vspace{0.2cm}

Enon\c{c}ons maintenant le deuxi\`eme r\'esultat principal de cette note.

\vspace{0.1cm}

\n {\bf Th\'eor\`eme II.1.} Soit $({\f g},B)$ une superalg\`ebre de Lie
quadratique qui est nilpotente ou r\'esoluble telle que
$[{\f g}_{\bar{1}},{\f g}_{\bar{1}}] \subset 
[{\f g}_{\bar{0}},{\f g}_{\bar{0}}].$
Alors, $\f g$ contient un id\'eal
gradu\'e totalement isotrope de dimension $E({dim{\f g}\over 2})$ qui est
maximal
parmi les sous-espaces vectoriels   totalement isotropes de $\f g$.
En plus, si la dimension de $\f g$ est paire alors $\f g$ est
isom\'etrique \`a une
$T^*-$extension de la superalg\`ebre de Lie quotient ${\f g}/{\f I}$. Si
la
dimension de $\f g$ est impaire, alors $\f g$ est isom\'etrique \`a un id\'eal
gradu\'e non
d\'eg\'en\'er\'e de codimension $1$ d'une $T^*-$extension de 
la superalg\`ebre de Lie quotient ${\f g}/{\f I}$.

\vspace{0.2cm}

\n {\bf Remarque.} Dans les d\'emonstrations du lemme II.1 et du th\'eor\`eme
II.1, nous
utilisons le th\'eor\`eme d'Engel dans le cas des superalg\`ebres de Lie
nilpotentes (o\`u il suffit de supposer que la caract\'eristique soit $\neq 2$)
 et le th\'eor\`eme de Lie dans le cas des superalg\`ebres de Lie r\'esolubles
$([4],[6])$. Rappelons le th\'eor\`eme de Lie dans le cas des 
superalg\`ebres de Lie
r\'esolubles: Toutes les repr\'esentations irr\'eductibles d'une superalg\`ebre
de Lie r\'esoluble ${\f g}= {\f g}_{\bar{0}}\oplus {\f g}_{\bar{1}}$ sur un
corps commutatif alg\'ebriquement clos de caract\'eristique z\'ero sont de
dimension $1$ si et seulement si
$[{\f g}_{\bar{1}},{\f g}_{\bar{1}}] \subset
                   [{\f g}_{\bar{0}},{\f g}_{\bar{0}}]$ 
([4], Prop.~5.2.4., p.82).
Ceci explique notre choix de la classe des superalg\`ebres de
Lie r\'esolubles consid\'er\'ee dans le lemme II.1 et dans le 
th\'eor\`eme II.1.

\vspace{0.2cm}

On termine cette note par les deux questions ouvertes suivantes:

1) Peut-on g\'en\'eraliser les r\'esultats du th\'eor\`eme II.1 \`a toutes les
superalg\`ebres de Lie r\'esolubles
${\f g}= {\f g}_{\bar{0}}\oplus {\f g}_{\bar{1}}$
sur un corps commutatif alg\'ebriquement clos de caract\'eristique z\'ero?

2) Dans $[5]$, il est montr\'e que toute alg\`ebre de Lie quadratique est
une double extension. Par cons\'equence, toute $T^*$-extension d'une alg\`ebre
de Lie est une double extension. Avec la g\'en\'eralisation de la notion de la
double extension ([1]) et la g\'en\'eralisation de la $T^*-$ extension
(dans cette note) aux superalg\`ebres de Lie, nous nous posons la question
naturelle suivante: Si $T^*_\o {\f g}$ est une $T^*-$extension d'une
superalg\`ebre de Lie $\f g$, $T^*_\o {\f g}$ est-elle une double
extension?

\vspace{0.1cm}

Les d\'emonstrations des r\'esultats annonc\'es para\^{\i}tront dans un article
ult\'erieur.

\vspace{0.4cm}

\n {\bf Remerciements}: I.B., S.B. et M.B. remercient les D\'epartements de 
 Math\'ematiques 
 de l'Universit\'e de Vigo et de Metz et le Graduiertenkolleg `Partielle 
 Differentialgleichungen' de l'Universit\'e de Freiburg pour des s\'ejours 
de recherche
 pendant lesquels ce travail a \'et\'e con\c{c}u.

\subsection*{R\'ef\'erences}

\begin{itemize}

 \item[[1]]  H. Benamor et S. Benayadi, Double extension of quadratic
            Lie superalgebras, {\sl Communications in Algebra}, 27 (1), 1999,
            p. 67-88.

\item[[2]]  S. Benayadi, Quadratic Lie superalgebras with the completely 
          reducible
          action of even part on the odd part, {\sl Journal of Algebra.}, 223,
           2000, p.344-366.

\item[[3]]  M. Bordemann, Nondegenerate invariant bilinear forms on
                nonassociative algebras, {\sl Acta Math. Univ. Comenianae,}
               Vol.~LXVI(1), 1997, p.~151-201.

\item [[4]]  V. Kac,   Lie superalgebras, {\sl Adv. Math}, 26, 1977, p. 8-96.

\item [[5]]  A. M\'edina  et  Ph. Revoy,   Alg\`ebres de Lie et produit
                  scalaire invariant, {\sl Ann. scient. Ec. Norm. Sup.}
                  4\`eme s\'erie, t.18, 1985, p.~553 - 561.

\item [[6]]  M. Scheunert, {\sl Theory of Lie superalgebras. An introduction,}
              Lecture Notes in Math. 716, 1979.

\end{itemize}

\end{document}